*Research Article*

# A Hybrid Natural Transform Homotopy Perturbation Method for Solving Fractional Partial Differential Equations

**Shehu Maitama**

*Department of Mathematics, Faculty of Science, Northwest University, Kano, Nigeria*

Correspondence should be addressed to Shehu Maitama; smusman12@sci.just.edu.jo

A hybrid analytical method for solving linear and nonlinear fractional partial differential equations is presented. The proposed analytical approach is an elegant combination of the Natural Transform Method (NTM) and a well-known method, Homotopy Perturbation Method (HPM). In this analytical method, the fractional derivative is computed in Caputo sense and the nonlinear term is calculated using He's polynomial. The proposed analytical method reduces the computational size and avoids round-off errors. Exact solution of linear and nonlinear fractional partial differential equations is successfully obtained using the analytical method.

## 1. Introduction

The concept of fractional calculus which deals with derivatives and integrals of arbitrary orders [1] plays a significant role in many areas of physical science and engineering. Recently, there is a rapid development in the concept of fractional calculus and its applications [2–5]. The linear and nonlinear fractional differential equations are used to model significant problems in fluid mechanics, acoustic, electromagnetism, signal processing, analytical chemistry, biology, and many other areas of physical science and engineering [6]. In recent years, many analytical and numerical methods have been used to solve linear and nonlinear fractional differential equations such as Adomian Decomposition Method [7–10], Homotopy Analysis Method [11], Laplace Decomposition Method [12], Homotopy Perturbation Method [13, 14], and Yang Laplace transform [15]. Moreover, local fractional variational iteration method [16, 17], Cylindrical-Coordinate method [18], modified Laplace Decomposition Method [12, 19], and fractional complex transform method [20] are also applied to linear and nonlinear fractional partial differential equations. The asymptotic behavior of solutions of a weighted Cauchy-type nonlinear fractional partial differential equations is studied in [21–23].

It is evident that most of the existing methods have so many deficiencies such as unnecessary linearization, discretization of variables, transformation, or taking some restrictive assumptions.

In this paper, an analytical method called a Hybrid Natural Transform Homotopy Perturbation Method for solving linear and nonlinear fractional partial differential equations is introduced. The proposed analytical method is applied directly without any linearization, transformation, discretization of variables, and so on. The analytical method gives a series solution which converges rapidly to an exact or approximate solution with elegant computational terms. The nonlinear terms are elegantly computed using He's polynomials [24–26]. Exact solution of linear and nonlinear fractional partial differential equation is successfully obtained using the new analytical method. Thus, the Hybrid Natural Transform Homotopy Perturbation Method is a powerful mathematical method for solving linear and nonlinear fractional partial differential equations and is a refinement of the existing methods.

## 2. Natural Transform

In this section, the basic definition and properties of the Natural Transform are presented.

In the year 2008, Z. H. Khan and W. A. Khan [27] introduced an integral transform called $N$-transform and it was recently renamed as the Natural Transform by Belgacem and Silambarasan [28, 29]. In fact, based on personal communications and Internet records, CASM in LUMS and IOSS Forman Christian College in Lahore lectures in 2004 and 2005 and 3rd ICM4F COMSATs University conference side lectures in Islamabad in 2006 indicate that Belgacem shared and discussed various aspects of this transform, with the various attending audiences. The Natural Transform is an integral transform which is similar to Laplace transform [30] and Sumudu integral transform [31, 32]. It converges to Laplace transform when $u = 1$ and Sumudu transform when $s = 1$. Belgacem and Silambarasan [29, 33] have proposed a detailed theory and applications of the Natural Transform. Recently, Natural Transform and Adomian Decomposition Method are successfully combined and obtained an exact solution of linear and nonlinear partial differential equations [34–39]. More details about the Natural Transform and its applications can be seen in Wikipedia note about the Natural Transform [40].

Table 1: List of some special natural transforms.

| Functional form | Natural transform form |
|---|---|
| 1 | $\dfrac{1}{s}$ |
| $t$ | $\dfrac{u}{s^2}$ |
| $e^{at}$ | $\dfrac{1}{s-au}$ |
| $\dfrac{t^{n-1}}{(n-1)!}, n = 1, 2, \ldots$ | $\dfrac{u^{n-1}}{s^n}$ |
| $\sin(t)$ | $\dfrac{u}{s^2+u^2}$ |

Some useful Natural Transforms in this paper are presented in Table 1.

*Definition 1.* The Natural Transform of the function $v(t) > 0$ and $v(t) = 0$ for $t < 0$ is defined over the set of functions:

$$A = \{v(t) : \exists M, \tau_1, \tau_2 > 0, \ |v(t)| < Me^{|t|/\tau_j}, \text{ if } t \in (-1)^j \times [0,\infty), \ j = 1, 2, \ldots\}, \quad (1)$$

by the following integral:

$$\mathbb{N}^+[v(t)] = V(s,u) = \frac{1}{u}\int_0^\infty e^{-st/u} v(t)\,dt; \quad s > 0, \ u > 0. \quad (2)$$

And the inverse Natural Transform of the function $v(t)$ is defined by

$$\mathbb{N}^{-1}[V(s,u)] = v(t) = \frac{1}{2\pi i}\int_{\gamma-i\infty}^{\gamma+i\infty} e^{st/u} V(s,u)\,ds, \quad (3)$$

where $s$ and $u$ are the Natural Transform variables [28, 29] and $\gamma$ is a real constant and the integral in (3) is taken along $s = \gamma$ in the complex plane $s = x + iy$.

Some properties of the Natural Transform Method are given below.

*Property 1.* Natural Transform of derivative: if $v^n(t)$ is the $n$th derivative of the function $v(t) \in A$ with respect to "$t$," then its Natural Transform is given by

$$\mathbb{N}^+[v^n(t)] = \frac{s^n}{u^n}V(s,u) - \sum_{k=0}^{n-1}\frac{s^{n-(k+1)}}{u^{n-k}}v^k(0). \quad (4)$$

When $n = 1, 2,$ and $3$, we have the following results:

$$\mathbb{N}^+[v'(t)] = \frac{s}{u}V(s,u) - \frac{1}{u}v(0),$$

$$\mathbb{N}^+[v''(t)] = \frac{s^2}{u^2}V(s,u) - \frac{s}{u^2}v(0) - \frac{1}{u}v'(0),$$

$$\mathbb{N}^+[v'''(t)] = \frac{s^3}{u^3}V(s,u) - \frac{s^2}{u^3}v(0) - \frac{s}{u^2}v'(0) - \frac{1}{u}v''(0). \quad (5)$$

*Property 2.* If $V(s,u)$ is the Natural Transform and $F(s)$ is the Laplace transform of the function $f(t) \in A$, then $\mathbb{N}^+[f(t)] = V(s,u) = (1/u)\int_0^\infty e^{-st/u} f(t)dt = (1/u)F(s/u)$.

*Property 3.* If $V(s,u)$ is the Natural Transform and $G(u)$ is the Sumudu transform of the function $v(t) \in A$, then $\mathbb{N}^+[v(t)] = V(s,u) = (1/s)\int_0^\infty e^{-t} v(ut/s)dt = (1/s)G(u/s)$.

*Property 4.* If $\mathbb{N}^+[v(t)] = V(s,u)$, then $\mathbb{N}^+[v(\beta t)] = (1/\beta)V(s/\beta, u)$.

*Property 5.* The Natural Transform is a linear operator; that is, if $\alpha$ and $\beta$ are nonzero constants, then

$$\mathbb{N}^+[\alpha f(t) \pm \beta g(t)] = \alpha\mathbb{N}^+[f(t)] \pm \beta\mathbb{N}^+[g(t)] = \alpha F^+(s,u) \pm \beta G^+(s,u). \quad (6)$$

Moreover, $F^+(s,u)$ and $G^+(s,u)$ are the Natural Transforms of functions $f(t)$ and $g(t)$, respectively.

## 3. Basic Definitions of Fractional Calculus

In this section, the basic definitions of fractional calculus are presented.

*Definition 2.* Function $f(t)$, $t > 0$, is said to be in the space $C^m_\alpha, m \in N \cup \{0\}$, if $f^{(m)} \in C_\alpha$.

*Definition 3.* A real function $f(t)$, $t > 0$, is said to be in the space $C_\alpha$, $\alpha \in \mathbb{R}$ if there exists a real number $p$ ($>\alpha$) such that $f(t) = t^p f_1(t)$, where $f_1(t) \in C[0,\infty)$. Clearly $C_\alpha \subset C_\beta$ if $\beta \leq \alpha$.

*Definition 4.* The left sided Riemann-Liouville fractional integral operator of order $\mu > 0$, of a function $f(t) \in C_\alpha$, and $\alpha \geq -1$ is defined as [19, 41]

$$I^\mu f(t) = \begin{cases} \dfrac{1}{\Gamma(\mu)} \int_0^t (t-\tau)^{\mu-1} f(\tau) d\tau, & \mu > 0, \ t > 0, \\ f(t), & \mu = 0, \end{cases} \quad (7)$$

where $\Gamma(\cdot)$ is the well-known Gamma function.

*Definition 5.* The left sided Caputo fractional derivative $f$, $f \in C_1^m, m \in \mathbb{N} \cup \{0\}$, is defined as [1, 5]

$$D_*^\mu f(t) = \begin{cases} I^{m-\mu}\left[\dfrac{\partial^m f(t)}{\partial t^m}\right], & m-1 < \mu < m, \ m \in \mathbb{N}, \\ \dfrac{\partial^m f(t)}{\partial t^m}, & \mu = m. \end{cases} \quad (8)$$

Note that [1, 5]

(i) $I_t^\mu f(t) = (1/\Gamma(\mu)) \int_0^t (f(t)/(t-s)^{1-\mu}) dt, \mu > 0, t > 0$,

(ii) $D_*^\mu f(x,t) = I_t^{m-\mu}[\partial^m f(t)/\partial t^m], m-1 < \mu \leq m$.

*Definition 6.* The Natural Transform of the Caputo fractional derivative is defined as

$$\mathbb{N}^+\left[D_t^{n\alpha} v(t)\right] = \frac{s^{n\alpha}}{u^{n\alpha}} V(s,u) - \sum_{k=0}^{n-1} \frac{s^{n\alpha-(k+1)}}{u^{n\alpha-k}} v^{(k)}(0+), \quad (9)$$

$((n-1)/n < \alpha \leq 1)$.

*Definition 7.* The series expansion defines a one-parameter Mittag-Leffler function as [1]

$$E_\alpha(z) = \sum_{k=0}^\infty \frac{z^k}{\Gamma(\alpha k + 1)}, \quad \alpha > 0, \ z \in \mathbb{C}. \quad (10)$$

## 4. Analysis of the Method

In this section, the basic idea of the Hybrid Natural Transform Homotopy Perturbation Method is clearly illustrated by the standard nonlinear fractional partial differential equation of the form

$$D_t^{n\alpha} v(x,t) + M(v(x,t)) + F(v(x,t)) = g(x,t), \quad (11)$$

subject to the initial condition

$$v(x,0) = f(x), \quad (12)$$

where $F(v(x,t))$ represents the nonlinear terms, $D_t^{n\alpha} = \partial^{n\alpha}/\partial t^{n\alpha}$ is the Caputo fractional derivative of function $v(t)$, $M(v(x,t))$ is the linear differential operator, and $g(x,t)$ is a source term.

Applying the Natural Transform to (11) subject to the given initial condition we get

$$V(x,s,u) = \frac{1}{s} f(x) + \frac{u^{n\alpha}}{s^{n\alpha}} \mathbb{N}^+ [g(x,t)] \\ - \frac{u^{n\alpha}}{s^{n\alpha}} \mathbb{N}^+ [M(v(x,t)) + F(v(x,t))]. \quad (13)$$

Taking the inverse Natural Transform of (13), we get

$$v(x,t) = G(x,t) \\ - \mathbb{N}^{-1}\left[\frac{u^{n\alpha}}{s^{n\alpha}} \mathbb{N}^+ [M(v(x,t)) + F(v(x,t))]\right], \quad (14)$$

where $G(x,t)$ is a term arising from the source term and the prescribed initial condition.

Now we apply the Homotopy Perturbation Method:

$$v(x,t) = \sum_{n=0}^\infty p^n v_n(x,t). \quad (15)$$

The nonlinear term $F(v(x,t))$ is decomposed as

$$F(v(x,t)) = \sum_{n=0}^\infty p^n H_n(v), \quad (16)$$

where $H_n(v)$ is He's polynomial which is computed using the following formula:

$$H_n(v_1, v_2, \ldots, v_n) = \frac{1}{n!} \frac{\partial^n}{\partial p^n}\left[F\left(\sum_{j=0}^n p^j v_j\right)\right]_{p=0}, \quad (17)$$

$$n = 0, 1, 2, \ldots.$$

Substituting (16) and (15) into (14), we get

$$\sum_{n=0}^\infty p^n v_n(x,t) = G(x,t) - p\left(\mathbb{N}^{-1}\left[\frac{u^{n\alpha}}{s^{n\alpha}}\right. \\ \left.\cdot \mathbb{N}^+\left[\sum_{n=0}^\infty p^n M(v(x,t)) + \sum_{n=0}^\infty p^n H_n(v)\right]\right]\right). \quad (18)$$

Using the coefficient of like powers of $p$ in (18), we obtained the following approximations:

$$p^0 : v_0(x,t) = G(x,t),$$

$$p^1 : v_1(x,t) \\ = -\mathbb{N}^{-1}\left[\frac{u^{n\alpha}}{s^{n\alpha}} \mathbb{N}^+ [M(v(x,t)) + H_0(v)]\right],$$

$$p^2 : v_2(x, t)$$
$$= -\mathbb{N}^{-1}\left[\frac{u^{n\alpha}}{s^{n\alpha}}\mathbb{N}^+\left[M(v(x,t)) + H_1(v)\right]\right],$$

$$p^3 : v_3(x, t)$$
$$= -\mathbb{N}^{-1}\left[\frac{u^{n\alpha}}{s^{n\alpha}}\mathbb{N}^+\left[M(v(x,t)) + H_2(v)\right]\right],$$

$$\vdots \quad (19)$$

and so on.

Hence, the series solution of (11)-(12) is given by

$$v(x, t) = \lim_{N \to \infty} \sum_{n=0}^{N} v_n(x, t). \quad (20)$$

## 5. Applications

In this section, the application of the Hybrid Natural Transform Homotopy Perturbation Method to linear and nonlinear fractional partial differential equations is clearly demonstrated to show its simplicity and high accuracy.

*Example 8.* Consider the following fractional diffusion equation of the form

$$D_t^\alpha v + v_{xx} + v_{yy} + v_{zz} = 0,$$
$$-\infty < x, y, z < \infty, \; t > 0, \quad (21)$$

subject to the initial condition

$$v(x, y, z, 0) = e^{x+y+z}, \quad \alpha \in (0, 1). \quad (22)$$

Applying the Natural Transform to (21) subject to the given initial condition, we get

$$V(x, y, z, s, u) = \frac{e^{x+y+z}}{s} + \frac{u^\alpha}{s^\alpha}\mathbb{N}^+\left[v_{xx} + v_{yy} + v_{zz}\right]. \quad (23)$$

Taking the inverse Natural Transform of (23), we get

$$v(x, y, z, t) = e^{x+y+z}$$
$$+ \mathbb{N}^{-1}\left[\frac{u^\alpha}{s^\alpha}\mathbb{N}^+\left[v_{xx} + v_{yy} + v_{zz}\right]\right]. \quad (24)$$

Now we apply Homotopy Perturbation Method:

$$v(x, y, z, t) = \sum_{n=0}^{\infty} p^n v_n(x, y, z, t). \quad (25)$$

Then (24) will become

$$\sum_{n=0}^{\infty} p^n v_n(x, y, z, t) = e^{x+y+z} - p\left(\mathbb{N}^{-1}\left[\frac{u^\alpha}{s^\alpha}\right.\right.$$
$$\left.\left.\cdot \mathbb{N}^+\left[\sum_{n=0}^{\infty} p^n v_{nxx} + \sum_{n=0}^{\infty} p^n v_{nyy} + \sum_{n=0}^{\infty} p^n v_{nzz}\right]\right]\right). \quad (26)$$

Using the coefficients of like powers of $p$ in (26), we obtained the following approximations:

$$p^0 : v_0(x, y, z, t) = e^{x+y+z},$$

$$p^1 : v_1(x, y, z, t)$$
$$= -\mathbb{N}^{-1}\left[\frac{u^\alpha}{s^\alpha}\mathbb{N}^+\left[v_{0xx} + v_{0yy} + v_{0zz}\right]\right]$$
$$= -\frac{e^{x+y+z}t^\alpha}{\Gamma(\alpha+1)},$$

$$p^2 : v_2(x, y, z, t)$$
$$= -\mathbb{N}^{-1}\left[\frac{u^\alpha}{s^\alpha}\mathbb{N}^+\left[v_{1xx} + v_{1yy} + v_{1zz}\right]\right] \quad (27)$$
$$= \frac{e^{x+y+z}t^{2\alpha}}{\Gamma(2\alpha+1)},$$

$$p^3 : v_3(x, y, z, t)$$
$$= -\mathbb{N}^{-1}\left[\frac{u^\alpha}{s^\alpha}\mathbb{N}^+\left[v_{2xx} + v_{2yy} + v_{2zz}\right]\right]$$
$$= -\frac{e^{x+y+z}t^{3\alpha}}{\Gamma(3\alpha+1)},$$

$$\vdots$$

and so on.

Then, the series solution of (21)-(22) is given by

$$v(x, y, z, t) = \lim_{N \to \infty} \sum_{n=0}^{N} v_n(x, y, z, t) = v_0(x, y, z, t)$$
$$+ v_1(x, y, z, t) + v_2(x, y, z, t) + v_3(x, y, z, t) + \cdots$$
$$= e^{x+y+z}\left(1 - \frac{t^\alpha}{\Gamma(\alpha+1)} + \frac{t^{2\alpha}}{\Gamma(2\alpha+1)}\right. \quad (28)$$
$$\left. - \frac{t^{3\alpha}}{\Gamma(3\alpha+1)} + \cdots\right) = e^{x+y+z}\left(1\right.$$
$$\left. + \sum_{m=1}^{\infty}\frac{(-t^\alpha)^m}{\Gamma(m\alpha+1)}\right) = e^{x+y+z}E_\alpha(-t^\alpha).$$

When $\alpha = 1$, the following result is obtained:

$$v(x, y, z, t) = \lim_{N \to \infty} \sum_{n=0}^{N} v_n(x, y, z, t)$$
$$= v_0(x, y, z, t) + v_1(x, y, z, t)$$
$$+ v_2(x, y, z, t) + v_3(x, y, z, t) + \cdots \quad (29)$$
$$= e^{x+y+z}\left(1 - \frac{t}{1!} + \frac{t^2}{2!} - \frac{t^3}{3!} + \cdots\right)$$
$$= e^{x+y+z-t},$$

which is the exact solution of (21)-(22), when $\alpha = 1$.

*Example 9.* Consider the following coupled system of nonlinear fractional partial differential equations of the form

$$D_t^\alpha v - v_{xx} - 2vv_x + v_x w_x = 0,$$
$$D_t^\alpha w - w_{xx} - 2ww_x + v_x w_x = 0,$$
(30)

subject to the initial conditions

$$v(x, 0) = \sin(x),$$
$$w(x, 0) = \sin(x).$$
(31)

Applying the Natural Transform to (30) subject to given initial conditions, we get

$$V(x, s, u) = \frac{1}{s}\sin(x) + \frac{u^\alpha}{s^\alpha}\mathbb{N}^+\left[v_{xx} + 2vv_x - v_x w_x\right],$$

$$W(x, s, u) = \frac{1}{s}\sin(x) + \frac{u^\alpha}{s^\alpha}\mathbb{N}^+\left[w_{xx} + 2ww_x - v_x w_x\right].$$
(32)

Taking the inverse Natural Transform of (32), we have

$$v(x, t) = \sin(x) + \mathbb{N}^{-1}\left[\frac{u^\alpha}{s^\alpha}\mathbb{N}^+\left[v_{xx} + 2vv_x - v_x w_x\right]\right],$$

$$w(x, t) = \sin(x) + \mathbb{N}^{-1}\left[\frac{u^\alpha}{s^\alpha}\mathbb{N}^+\left[ww_{xx} + 2ww_x - v_x w_x\right]\right].$$
(33)

Now we apply the Homotopy Perturbation Method:

$$v(x, t) = \sum_{n=0}^{\infty} p^n v_n(x, t),$$

$$w(x, t) = \sum_{n=0}^{\infty} p^n w_n(x, t).$$
(34)

Then (33) will become

$$v(x, t) = \sin(x) + p\left(\mathbb{N}^{-1}\left[\frac{u^\alpha}{s^\alpha}\right.\right.$$
$$\left.\left.\cdot \mathbb{N}^+\left[\sum_{n=0}^{\infty} p^n v_{nxx} + 2\sum_{n=0}^{\infty} p^n H_n - \sum_{n=0}^{\infty} p^n H_n''\right]\right]\right),$$

$$w(x, t) = \sin(x) + p\left(\mathbb{N}^{-1}\left[\frac{u^\alpha}{s^\alpha}\right.\right.$$
$$\left.\left.\cdot \mathbb{N}^+\left[\sum_{n=0}^{\infty} p^n w_{nxx} + 2\sum_{n=0}^{\infty} p^n H_n' - \sum_{n=0}^{\infty} p^n H_n''\right]\right]\right),$$
(35)

where $H_n$, $H_n'$, and $H_n''$ are He's polynomials which represent the nonlinear terms $vv_x$, $ww_x$, and $v_x w_x$, respectively.

Some few components of He's polynomials of the nonlinear terms $vv_x$, $ww_x$, and $v_x w_x$ are given as follows:

$$H_0 = v_0 v_{0x},$$
$$H_1 = v_0 v_{1x} + v_1 v_{0x},$$
$$H_2 = v_0 v_{2x} + v_1 v_{1x} + v_2 v_{0x},$$
$$\vdots$$
$$H_0' = w_0 w_{0x},$$
$$H_1' = w_0 w_{1x} + w_1 w_{0x},$$
$$H_2' = w_0 w_{2x} + w_1 w_{1x} + w_2 w_{0x},$$
$$\vdots$$
$$H_0'' = v_{0x} w_{0x},$$
$$H_1'' = v_{0x} w_{1x} + v_{1x} w_{0x},$$
$$H_2'' = v_{2x} w_{0x} + v_{1x} w_{1x} + v_{0x} w_{2x},$$
$$\vdots$$
(36)

and so on.

Using the coefficients of the like powers of $p$ in (35), we obtained the following approximations:

$$p^0 : v_0(x, t) = \sin(x),$$

$$p^1 : v_1(x, t) = \mathbb{N}^{-1}\left[\frac{u^\alpha}{s^\alpha}\mathbb{N}^+\left[v_{0xx} + 2H_0 - H_0''\right]\right]$$
$$= -\frac{\sin(x) t^\alpha}{\Gamma(\alpha + 1)},$$

$$p^2 : v_2(x, t) = \mathbb{N}^{-1}\left[\frac{u^\alpha}{s^\alpha}\mathbb{N}^+\left[v_{1xx} + 2H_1 - H_1''\right]\right]$$
$$= \frac{\sin(x) t^{2\alpha}}{\Gamma(2\alpha + 1)},$$

$$p^3 : v_3(x, t) = \mathbb{N}^{-1}\left[\frac{u^\alpha}{s^\alpha}\mathbb{N}^+\left[v_{2xx} + 2H_2 - H_2''\right]\right]$$
$$= -\frac{\sin(x) t^{3\alpha}}{\Gamma(3\alpha + 1)},$$

$$\vdots$$
(37)

and so on.

Similarly,

$$p^0 : w_0(x, t) = \sin(x),$$

$$p^1 : w_1(x, t) = \mathbb{N}^{-1}\left[\frac{u^\alpha}{s^\alpha}\mathbb{N}^+\left[w_{0xx} + 2H_0' - H_0''\right]\right]$$
$$= -\frac{\sin(x) t^\alpha}{\Gamma(\alpha + 1)},$$

$$p^2 : w_2(x,t) = \mathbb{N}^{-1}\left[\frac{u^\alpha}{s^\alpha}\mathbb{N}^+\left[w_{1xx} + 2H_1' - H_1''\right]\right]$$
$$= \frac{\sin(x)\,t^{2\alpha}}{\Gamma(2\alpha+1)},$$
$$p^3 : w_3(x,t) = \mathbb{N}^{-1}\left[\frac{u^\alpha}{s^\alpha}\mathbb{N}^+\left[w_{2xx} + 2H_2' - H_2''\right]\right]$$
$$= -\frac{\sin(x)\,t^{3\alpha}}{\Gamma(3\alpha+1)},$$
$$\vdots \tag{38}$$

and so on.

Thus, the series solution of (30)-(31) is given by

$$v(x,t) = \lim_{N\to\infty}\sum_{n=0}^{N} v_n(x,t) = v_0(x,t) + v_1(x,t)$$
$$+ v_2(x,t) + \cdots = \sin(x)$$
$$\cdot\left(1 - \frac{t^\alpha}{\Gamma(\alpha+1)} + \frac{t^{2\alpha}}{\Gamma(2\alpha+1)} - \frac{t^{3\alpha}}{\Gamma(3\alpha+1)} + \cdots\right)$$
$$= \sin(x)\left(1 + \sum_{m=1}^{\infty}\frac{(-t^\alpha)^m}{\Gamma(m\alpha+1)}\right) = \sin(x)E_\alpha(-t^\alpha),$$
$$w(x,t) = \lim_{N\to\infty}\sum_{n=0}^{N} w_n(x,t) = w_0(x,t) + w_1(x,t)$$
$$+ w_2(x,t) + \cdots = \sin(x)$$
$$\cdot\left(1 - \frac{t^\alpha}{\Gamma(\alpha+1)} + \frac{t^{2\alpha}}{\Gamma(2\alpha+1)} - \frac{t^{3\alpha}}{\Gamma(3\alpha+1)} + \cdots\right)$$
$$= \sin(x)\left(1 + \sum_{m=1}^{\infty}\frac{(-t^\alpha)^m}{\Gamma(m\alpha+1)}\right) = \sin(x)E_\alpha(-t^\alpha). \tag{39}$$

Hence, the exact solutions of (30)-(31) are given by
$$v(x,t) = \sin(x)E_\alpha(-t^\alpha),$$
$$w(x,t) = \sin(x)E_\alpha(-t^\alpha). \tag{40}$$

When $\alpha = 1$, the following result is obtained:
$$v(x,t) = \lim_{N\to\infty}\sum_{n=0}^{N} v_n(x,t)$$
$$= v_0(x,t) + v_1(x,t) + v_2(x,t) + \cdots$$
$$= \sin(x)\left(1 - t + \frac{t^2}{2!} + \cdots\right) = e^{-t}\sin(x),$$
$$w(x,t) = \lim_{N\to\infty}\sum_{n=0}^{N} w_n(x,t)$$
$$= w_0(x,t) + w_1(x,t) + w_2(x,t) + \cdots$$
$$= \sin(x)\left(1 - t + \frac{t^2}{2!} + \cdots\right) = e^{-t}\sin(x), \tag{41}$$

which is the exact solution of (30)-(31), when $\alpha = 1$.

## 6. Conclusion

In this paper, Natural Transform Method (NTM) and Homotopy Perturbation Method (HPM) are successfully combined to form a robust analytical method called a Hybrid Natural Transform Homotopy Perturbation Method for solving linear and nonlinear fractional partial differential equations. The analytical method gives a series solution which converges rapidly to an exact or approximate solution with elegant computational terms. In this analytical method, the fractional derivative is computed in Caputo sense, while the nonlinear term is calculated using He's polynomials. The analytical procedure is applied successfully and obtained an exact solution of linear and nonlinear fractional partial differential equations. The simplicity and high accuracy of the analytical method are clearly illustrated. Thus, the Hybrid Natural Transform Homotopy Perturbation Method is a powerful analytical method for solving linear and nonlinear fractional partial differential equations.

## Competing Interests

The author declares that there is no conflict of interests regarding the publication of this article.

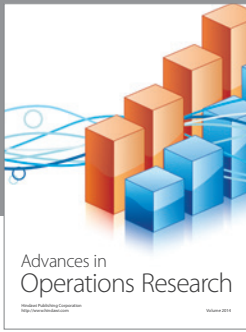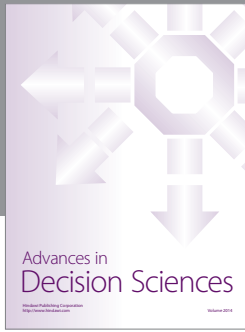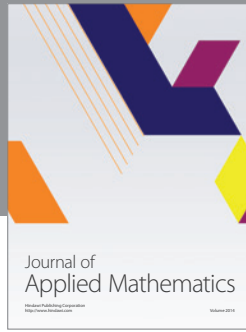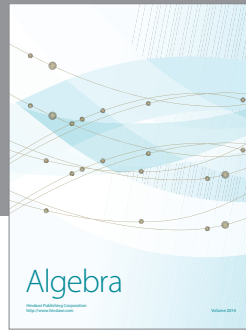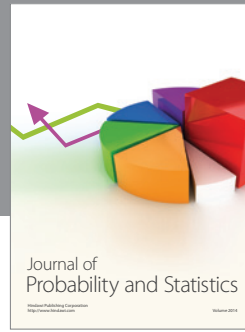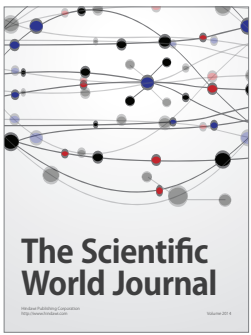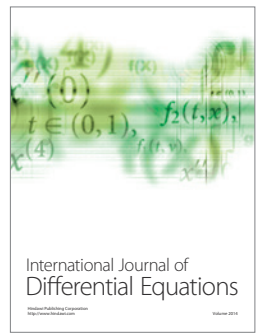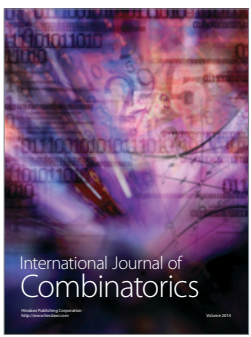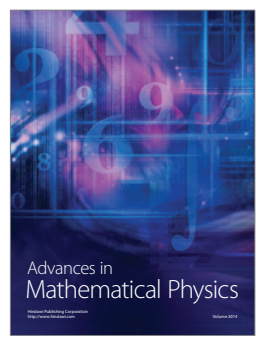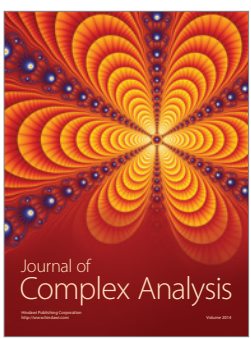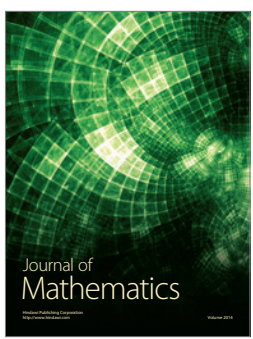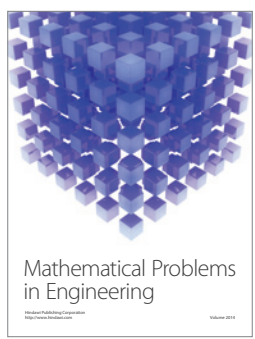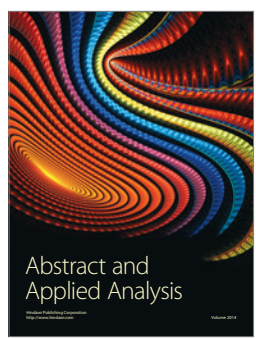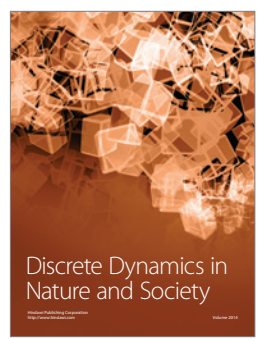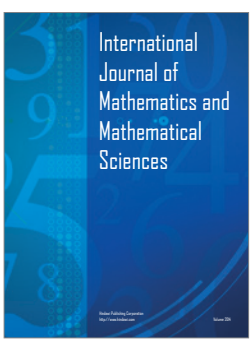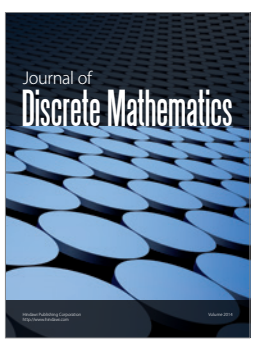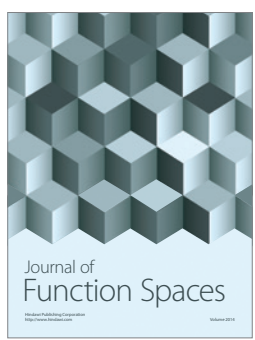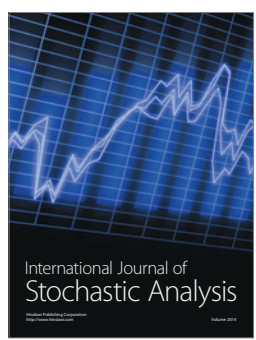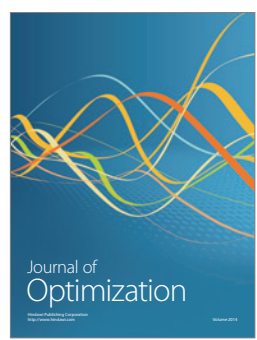

Submit your manuscripts at
http://www.hindawi.com